\documentclass[12pt, reqno]{amsart}
\usepackage{amsmath}
\usepackage[mathscr]{eucal}
\usepackage{amssymb}
\usepackage{latexsym}
\usepackage{amsthm}
\usepackage{graphicx}
\usepackage{bm}

\theoremstyle{plain}

\newtheorem{thm}{Theorem}[section]
\newtheorem{cor}{Corollary}[section]
\newtheorem{lem}{Lemma}[section]

\theoremstyle{definition}
\newtheorem{rem}{Remark}[section]
\newtheorem{dfn}{Definition}[section]
\newtheorem{prf}{Proof}[section]

\title{On almost convergence on locally compact abelian groups}
\author{Ryoichi Kunisada}

\address{Faculty of Liberal Arts, Tsuru University, Tsuru-shi, Yamanashi-ken 402-8555, Japan}

\email{tk-waseda@ruri.waseda.jp}

\thanks{}



\date{}

\begin{document}
\maketitle

\begin{abstract}
We study a summability method called almost convergence for bounded measurable functions defined on a locally compact abelian group. We define almost convergence using topologically invariant means and exhibit two different kinds of necessary and sufficient conditions, one is analytic and the other is functional analytic, for a given function to be almost convergent. As an application, we show complex Tauberian theorems for almost convergence on the integers and the real numbers. These results are closely related to some of the classical Tauberian theorems like the Wiener-Ikehara and Katznelson-Tzafriri theorems.
\end{abstract}

\bigskip

\section{Introduction}
For a locally compact abelian group $G$, let $L^1(G)$ be the group algebra of $G$ and $L^{\infty}(G)$ be the set of all essentially bounded measurable functions on $G$. Let $C_{bu}(G)$ be the set of all bounded, uniformly continuous functions on $G$. The uniformity $\mathcal{U}$ of $G$ is the set of all subsetes of $G^2$ given by 
\[\{(x, y) \in G^2 : x-y \in U\}, \]
where $U$ is a neighborhood of $0$. Note that both spaces $L^{\infty}(G)$ and $C_{bu}(G)$ are Banach spaces with respect to the supremum norm $\|\cdot \|_{\infty}$. Here we consider complex-valued functions in general and we denote by $L^{\infty}_R(G)$ the space of real-valued essentially bounded measurable functions.

A general element of $L^1(G)$ is denotedy by the symbol $f$(we also use $g, h$ if necessary) and that of $L^{\infty}(G)$ is denoted by $\psi$. For each $s \in \mathbb{R}$, we use the symbols $f_s(x) := f(x+s)$ and $\psi_s(x) := \psi(x+s)$, the translates of $f$ and $\psi$ by $s$, respectively.

Let $L^{\infty}(G)^*$ and $C_{bu}(G)^*$ be the dual spaces of $L^{\infty}(G)$ and $C_{bu}(G)$, respectively.
An element $\varphi$ of $L^{\infty}(G)^* \ (C_{bu}(G)^*)$ is said to be a {\it mean} on $L^{\infty}(G) \ (C_{bu}(G))$ if it satisfies
 
\smallskip
\noindent
$(1)$ $\varphi \ge 0 \ (\text{i.e.}, \varphi(\psi) \ge 0 \ \text{for every positive} \ \psi \in L^{\infty}(G) \ (C_{bu}(G)))$; \\
$(2)$ $\varphi(\bm{1}) = 1$, where $1$ is the constant function taking the value $1$ everywhere.

\smallskip
\noindent
Further, $\varphi$ is called an {\it invariant mean} on $L^{\infty}(G) \ (C_{bu}(G))$ if it is a mean such that 

\smallskip
\noindent
$(3)$ $\varphi(\psi_s) = \varphi(\psi)$ for every $\psi \in L^{\infty}(G) \ (C_{bu}(G))$ and $s \in G$.

\smallskip
\noindent
Let us denote by $\mathcal{I}(G) \ (\mathcal{I}_0(G))$ the set of all invariant means on $L^{\infty}(G) \ (C_{bu}(G))$. 

Now, we introduce another class of means on $L^{\infty}(G)$ which is more important for our purpose. For $f \in L^1(G)$ and $\psi \in L^{\infty}(G)$, their convolution $f * \psi$ is defined by
\[f * \psi(x) = \int_G f(t)\psi(x-t)dm(t), \quad x \in G,  \]
where $m$ is the Haar measure of $G$. Let $P(G)$ be the set of positive elements $f$ in $L^1(G)$ such that $\int_G f(x)dm(x) = 1$. A mean $\varphi$ on $L^{\infty}(G)$ is said to be a {\it topologically invariant mean} if it satisfies the following condition (\cite{Hul}):

\smallskip
\noindent
$(4)$ $\varphi(f * \psi) = \varphi(\psi) \ \forall \psi \in L^{\infty}(G) \ \forall f \in P(G)$.

\smallskip
Let us denote by $\mathcal{T}(G)$ the set of all topologically invariant means on $L^{\infty}(G)$. Note that a topologically invariant mean is an invariant mean. In fact, if $\varphi$ is topologically invariant, then 
\[\varphi(\psi_s) = \varphi(f * \psi_s) = \varphi(f_s * \psi) = \varphi(\psi) \]
for each $\psi \in L^{\infty}(\mathbb{R})$ and $s \in \mathbb{R}$. Conversely, For discrete groups $G$, it is easy to show that $\mathcal{I}(G) = \mathcal{T}(G)$ holds true. For every locally compact abelian group $G$, $\mathcal{T}(G) \not= \emptyset$ is valid, and if $G$ is continuous, we have $\mathcal{T}(G) \subsetneq \mathcal{I}(G)$. For compact abelian groups $G$, $\mathcal{T}(G)$ is the singleton consisting of the normalized Haar measure of $G$. We refer the reader to \cite{Green}, \cite{Pat} for more detailed exposition on these results.

The main objective of this paper is a summability method concerning (topologically) invariant means. Let $l_{\infty}$ be the set of all bounded functions on the nonnegative integers $\mathbb{Z}_+ := \{n \in \mathbb{Z} : n \ge 0\}$. Lorentz (1948) defined a summability method on $l_{\infty}$ called {\it almost convergence} using Banach limits (\cite{Lor}). Recall that an element $\varphi$ of $l_{\infty}^*$, the dual space of $l_{\infty}$, is called a {\it Banach limit} if the following conditions are satisfied: 

\smallskip
\noindent
$(1)$ $\varphi \ge 0$; \\
$(2)$ $\varphi(\bm{1}) = 1$, where $1$ is the constant function taking the value $1$ everywhere. \\
$(3)$ $\varphi(\psi_n) = \varphi(\psi)$ for every $n \in \mathbb{Z}_+$ and $\psi \in l_{\infty}$, 

\smallskip
\noindent
namely, $\varphi$ is a right translation invariant mean on $l_{\infty}$.  Let $\mathcal{B}$ be the set of all Banach limits. Then, almost convergence for sequences is defined as follows.

 
\begin{dfn}[Lorentz, 1948]
$\psi \in l_{\infty}$ is said to be almost convergent to a (complex) number $\alpha$ if $\varphi(\psi) = \alpha$ for every $\varphi \in \mathcal{B}$.
\end{dfn}

It is difficult to know from this abstract definition whether a given sequence of numbers is almost convergent or not. However, Lorentz proved an analytic condition for almost convergence as follows (his delivation was rather complicated and Sucheston \cite{Suc} gave a more simple proof later). 
\begin{thm}[Lorentz, 1948]
$\psi \in l_{\infty}$ is almost convergent to $\alpha$ if and only if 
\[\lim_{k \to \infty} \frac{1}{k} \sum_{i=0}^{k-1} \psi(n+i) = \alpha \]
uniformly in $n \in \mathbb{Z}_+$.
\end{thm}

This assertion can be proved by the Hahn-Banach theorem and is a basic way to check whether a given sequence actually almost converges.

Note that Banach limits can be regarded as a special kind of invariant means on $L^{\infty}(\mathbb{Z})$, where $\mathbb{Z}$ is the additive group of integers. Hence, it is natural to consider that one can define the notion of almost convergence on an arbitrary locally compact abelian group or its subsemigroups (note that abelian groups $G$ are amenable, that is, there exist invariant means on $L^{\infty}(G)$). Specifically, one can define formally almost convergence of functions in $L^{\infty}(G)$ by the following :
\begin{dfn}
Let $G$ be a locally compact abelian group. We say that $\psi \in L^{\infty}(G)$ is almost convergent to a complex number $\alpha$ if and only if
\[\varphi(\psi) = \alpha \]
holds for every $\varphi \in \mathcal{T}(G)$. In this case, we write as $\psi \xrightarrow{ac} \alpha$.
\end{dfn}

The reason for adopting topological invariant means instead of (seemingly more natural) invariant means is that an analogous result of Theorem $1.1$ is also valid for groups not necessarily discrete (see \cite{Chow1}, \cite{Chow2}). 

The objective of this paper is to provide a new necessary and sufficient condition for a given bounded function on a locally compact abelian group to be almost convergent. This is obtained  through the thoery of harmonic analysis on locally compact abelian groups, especially, the theory of spectral synthesis. Furthermore, applying this result, we obtain complex tauberian theorems for almost convergence on $\mathbb{Z}$ and $\mathbb{R}$, which are related to the famous results of Katznelson-Tzafriri and Wiener-Ikehara.

The paper is organized as follows. In Section 2, we exhibit an analytic condition for almost convergence. This is a special case of the more general result of Chow \cite{Chow1} when the underling group is abelian. However, considering the importance of the result, we include this section for the sake of completeness and consistency of the paper. In Section 3, Abelian and Tauberian theorems concerning almost convergence are treated.

In Section 4, following Forelli \cite{Fore}, we develop spectral theory of $C_{bu}(G)$ and $C_{bu}(G)^*$, which contains a description of subspaces of $C_{bu}(G)^*$ defined via spectrum as the annihilator of a certain invariant subspace of $C_{bu}(G)$. This is somewhat a formal generalization of his relult to general locally compact abelian groups. In Section 5, as an application of a result of the previous section, we provide the annihilator of the subspace of $L^{\infty}(G)^*$ spanned by topologically invariant means on $L^{\infty}(G)$. In Section 6, using a result of Section 5, we obtain a new necessary and sufficient condition for almost convergence. Section 7 deals with almost convergence on positive subsemigroups of the special groups $\mathbb{Z}$ and $\mathbb{R}$. In Section 8, we deal with complex Tauberian theorems for almost convergence on $\mathbb{Z}$ and $\mathbb{R}$. 

\section{An analytic condition for almost convergence}
Let $G$ be a locally compact abelian group. Since $G$ is an amenable group, there exists a {\it summing net for $G$}, namely, a net $\{K_{\delta}\}_{\delta \in \Delta}$ of nonnull, compact subsets of $G$ satisfying the following properties (see \cite{Green}, \cite{Pat}): 

\smallskip
\noindent
$(1)$ $K_{\delta} \subseteq K_{\delta^{\prime}}$ if $\delta \le \delta^{\prime}$ \\
$(2)$ $G = \bigcup_{\delta \in \Delta} K_{\delta}$ \\
$(3)$ $\lim_{\delta} \frac{m(K_{\delta} \bigtriangleup K_{\delta, s})}{m(K_{\delta})} = 0$ uniformly in $s$ on a compact subset of $G$.

\smallskip
We fix one summing net for $G$ and define a sublinear functional on $L^{\infty}(G)$ as
\[\overline{p}(\psi) := \limsup_{\delta} \sup_{x \in G} \frac{1}{m(K_{\delta})} \int_{K_{\delta}} \psi_x(t)dt, \]
where $\psi \in L^{\infty}(G)$. We also introduce the functional $\underline{p}(\psi) := -\overline{p}(-\psi)$, which can be expressed as
\[\underline{p}(\psi) := \liminf_{\delta} \inf_{x \in G}  \frac{1}{m(K_{\delta})} \int_{K_{\delta}} \psi_x(t)dt. \]
Then, we have the following result:
\begin{thm}
A mean $\varphi$ on $C_{bu}(G)$ is an ivariant mean if and only if 
\[\varphi(\psi) \le \overline{p}(\psi) \]
holds for every $\psi \in C_{bu}(G)$.
\end{thm}

\begin{prf}
First, we show necessity. Let $\varphi$ be an invariant mean on $C_{bu}(G)$. Note that, since $\psi$ is in $C_{bu}(G)$, the mapping $G \ni t \mapsto \psi_t \in C_{bu}(G)$ is continuous and thus, it is Bochner integrable. Then, we have 
\[\varphi(\psi) = \frac{1}{m(K_{\delta})} \int_{K_{\delta}} \varphi(\psi_t)dt = \varphi\left(\frac{1}{m(K_{\delta})} \int_{K_{\delta}} \psi_t(x)dt\right) \le \sup_x \frac{1}{m(K_{\delta})} \int_{K_{\delta}} \psi_x(t)dt  \]
for every $K_{\delta}$ (see \cite{Yos}). Here we use the elementary fact that for any mean $\varphi$ on $C_{bu}(\mathbb{R})$ and $\psi$ in $C_{bu}(\mathbb{R})$, it holds that
\[\varphi(\psi) \le \sup_{x \in \mathbb{R}} \psi(x). \]
In fact, let $\alpha := \sup_{x \in \mathbb{R}} \psi(x)$ and then $\alpha - \psi \ge 0$, thus by the positivity of $\varphi$, we have
\begin{align}
\varphi(\alpha - \psi) \ge 0 &\Leftrightarrow \varphi(\alpha) \ge \varphi(\psi) \notag \\
&\Leftrightarrow \alpha \ge \varphi(\psi). \notag
\end{align}
Taking the limit superior over ${K_{\delta}}'s$, we obtain 
\[\varphi(\psi) \le \limsup_{\delta} \sup_x \frac{1}{m(K_{\delta})} \int_{K_{\delta}} \psi_x(t)dt. \]

Now, we show sufficiency. Let $\varphi$ satisfy the condition in the theorem. Then, for each $\psi$ in $C_{bu}(G)$, we have
\begin{align}
\varphi(\psi-\psi_s) &\le \limsup_{\delta} \sup_x \frac{1}{m(K_{\delta})} \int_{K_{\delta}} (\psi_x(t) - \psi_{x+s}(t))dt \notag \\
&= \limsup_{\delta} \sup_x \frac{1}{m(K_{\delta})} \left(\int_{K_{\delta}} \psi_x(t)dt - \int_{K_{\delta}}\psi_{x+s}(t)dt\right) \notag \\
&= \limsup_{\delta} \sup_x \frac{1}{m(K_{\delta})} \left(\int_{K_{\delta}} \psi_x(t)dt - \int_{K_{\delta, -s}}\psi_x(t)dt\right)  \notag \\
&\le \limsup_{\delta} \sup_x \frac{1}{m(K_{\delta})} \int_{K_{\delta} \bigtriangleup K_{\delta, -s}} |\psi_x(t)|dt \notag \\
&\le \limsup_{\delta} \frac{m(K_{\delta} \bigtriangleup K_{\delta, -s})}{m(K_{\delta})} \|\psi\|_{\infty} = 0. \notag
\end{align}
For the reverse inequality, we consider the relation 
\[\varphi(-\psi) \le \overline{p}(-\psi) \Leftrightarrow \varphi(\psi) \ge -\overline{p}(-\psi) =: \underline{p}(\psi). \]
Note that we can show $\underline{p}(\psi-\psi_s) = 0$ in the same argument as above and we obtain 
\[\varphi(\psi - \psi_s) \ge 0. \]
Hence, we conclude that $\varphi(\psi - \psi_s) = 0$, which means translation invariance of $\varphi$ and we complete the proof.

\end{prf}

We need the following lemma to obtain a similar condition for topologically invraiant means. An equivalent assertion was shown in \cite{Chow1} in a more general form.
\begin{lem}
For any $\psi \in L^{\infty}(G)$ and $f \in P(G)$, 
\[\overline{p}(\psi - f * \psi) = \underline{p}(\psi-f*\psi) = 0 \]
holds true.
\end{lem}

\begin{prf}
By direct computation, we have
\begin{align}
\overline{p}(\psi - f * \psi) &= \limsup_{\delta} \sup_{x \in G} \frac{1}{m(K_{\delta})} \int_{K_{\delta}} \{\psi_x(t) - (f * \psi)_x(t)\}dt \notag \\
&= \limsup_{\delta} \sup_{x \in G} \frac{1}{m(K_{\alpha})} \int_{K_{\delta}} \int_G \{\psi_x(t) - \psi_x(t-u)\}f(u)dudt \notag \\ 
&=  \limsup_{\delta} \sup_{x \in G} \int_G f(u)du \frac{1}{m(K_{\delta})} \int_{K_{\delta}} \{\psi_x(t) - \psi_x(t-u)\}dt. \notag
\end{align}
The integral 
\[\frac{1}{m(K_{\delta})} \int_{K_{\delta}} \{\psi_x(t) - \psi_x(t-u)\}dt \]
can be evaluated in two ways:
\begin{align}
\frac{1}{m(K_{\delta})}\int_{K_{\delta}} \{\psi_x(t) - \psi_x(t-u)\}dt &\le \left|\frac{1}{m(K_{\delta})}\int_{K_{\delta} \bigtriangleup K_{\delta, -u}} \psi_x(t)dt\right|   \notag \\
&\le \frac{m(K_{\delta} \bigtriangleup K_{\delta, -u})}{m(K_{\delta})} \|\psi\|_{\infty}. 
\end{align}
Also, we obtain 
\begin{equation}
\frac{1}{m(K_{\delta})}\int_{K_{\delta}} \{\psi_x(t) - \psi_x(t-u)\}dt \le \frac{1}{m(K_{\delta})} 2\|\psi\|_{\infty} m(K_{\delta}) = 2\|\psi\|_{\infty}. 
\end{equation}
Let $\varepsilon > 0$ be given. Take a compact subset $C_{\varepsilon}$ such that $\int_{G \setminus C_{\varepsilon}} |f(t)|dt < \varepsilon$. Then, we can choose $\gamma$ such that $\frac{m(K_{\delta} \bigtriangleup K_{\delta, -u})}{m(K_{\delta})} \le \varepsilon$ for every $ \delta \ge \gamma$ and $u \in C_{\varepsilon}$. Then, by $(1)$ and $(2)$, we have
\begin{align}
&\limsup_{\delta} \sup_{x \in G} \int_G f(u)du \frac{1}{m(K_{\delta})} \int_{K_{\delta}} \{\psi_x(t) - \psi_x(t-u)\}dt \notag \\
&= \limsup_{\delta} \sup_{x \in G} \bigg\{\int_{C_{\varepsilon}} f(u)du \frac{1}{m(K_{\delta})} \int_{K_{\delta}} \{\psi_x(t) - \psi_x(t-u)\}dt \notag \\
&+ \int_{G \setminus C_{\varepsilon}} f(u)du \frac{1}{m(K_{\delta})} \int_{K_{\delta}} \{\psi_x(t) - \psi_x(t-u)\}dt\bigg\} \notag \\
&\le \limsup_{\delta} \sup_{x \in G} \int_{C_{\varepsilon}} |f(u)| \frac{m(K_{\delta} \bigtriangleup K_{\delta, -u})}{m(K_{\delta})} \|\psi\|_{\infty} du + 2\|\psi\|_{\infty} \int_{G \setminus C_{\varepsilon}} |f(u)|du  \notag \\
&\le \varepsilon \|\psi\|_{\infty} \int_{C_{\varepsilon}} |f(u)|du + 2\varepsilon\|\psi\|_{\infty} \le 3\varepsilon \|\psi\|_{\infty}. \notag
\end{align}
Since $\varepsilon > 0$ can be arbitrary, we obtain $\overline{p}(\psi - f*\psi) \le 0$. The equation $\underline{p}(\psi - f * \psi) \ge 0$ can be proved similarly.  Since $\underline{p}(\psi) \le \overline{p}(\psi)$ holds for each $\psi \in L^{\infty}(G)$, we have $\overline{p}(\psi-f*\psi) = \underline{p}(\psi-f*\psi) = 0$. We complete the proof.
\end{prf}

Combining Theorem 2.1 and Lemma 2.1, we obtain the following result.
\begin{thm}
A mean $\varphi$ on $L^{\infty}(G)$ is a topologically ivariant mean if and only if 
\[\varphi(\psi) \le \overline{p}(\psi) \]
holds for every $\psi \in L^{\infty}(G)$.
\end{thm}

\begin{prf}
First, we show necessity. Assume that $\varphi$ is in $\mathcal{T}(G)$. Then, for any $\psi \in L^{\infty}(G)$ and $f \in P(G)$, we have 
\[\varphi(f * \psi) = \varphi(\psi). \]
Note that, as stated before, $\varphi$ is an invariant mean on $C_{bu}(G)$. Since $f * \psi \in C_{bu}(G)$, by Theorem $2.1$, we have
\[\varphi(f * \psi) \le \overline{p}(f * \psi). \]
Also, by Lemma $2.1$, $\overline{p}(f * \psi) = \overline{p}(\psi)$ holds and we obtain
\[\varphi(\psi) = \varphi(f * \psi) \le \overline{p}(f * \psi) = \overline{p}(\psi). \]

Now, we show sufficiency. Assume that $\varphi(\psi) \le \overline{p}(\psi)$ holds for each $\psi \in L^{\infty}(G)$. Then, by Lemma $2.1$, for any $f \in P(G)$, we have
\[0 = \underline{p}(\psi-f*\psi) \le \varphi(\psi - f * \psi) \le \overline{p}(\psi - f * \psi) = 0. \]
Hence, we obtain $\varphi(\psi - f * \psi) = 0$ and $\varphi$ is topologically invariant.
\end{prf}

\begin{thm}
Let $\psi \in L^{\infty}_{\mathbb{R}}(G)$. Then, $\psi$ is almost convergent to $\alpha$ if and only if 
\[\lim_{\delta} \frac{1}{m(K_{\delta})} \int_{K_{\delta}} \psi(x+s)dx = \alpha \]
uniformly in $s \in G$. 

\end{thm}

\begin{prf}
Note that, for any $\varphi \in \mathcal{T}(G)$ and $\psi \in L^{\infty}(G)$, we have
\[\underline{p}(\psi) \le \varphi(\psi) \le \overline{p}(\psi). \]
Conversely, for any real number $\alpha$ with $\underline{p}(\psi) \le \alpha \le \overline{p}(\psi)$, there exists some $\varphi \in \mathcal{T}(G)$ such that $\varphi(\psi) = \alpha$. In fact, we define $\varphi_0$ on $\mathbb{R}\psi = \{c\psi : c \in \mathbb{R}\}$ by $\varphi_0(c\psi) = c\alpha$ and then can extend it to whole $L^{\infty}(G)$ such that $\varphi(\psi) \le \overline{p}(\psi)$ holds for every $\psi \in L^{\infty}(G)$ by the Hahn-Banach theorem. This extended $\varphi$ is an topologically invariant mean by Theorem 2.2. 

Hence, we have shown that $\varphi(\psi) = \alpha$ for every $\varphi \in \mathcal{T}(G)$ if and only if $\underline{p}(\psi) = \overline{p}(\psi) = \alpha$. Now, we show that this is equivalent to the condition given in the theorem. First, necessity is clear.
Recall that the functionals $\overline{p}$ and $\underline{p}$ are expressed as folllows:
\[\overline{p}(\phi) = \limsup_{\delta} \sup_{s \in G} \frac{1}{m(K_{\delta})} \int_{K_{\delta}} \psi(x+s)dx, \]
\[\underline{p}(\psi) = \liminf_{\delta} \inf_{s \in G} \frac{1}{m(K_{\delta})} \int_{K_{\delta}} \psi(x+s)dx. \]
Let $\varepsilon > 0$ be any positive number. Then, there exist indicies $\delta_0$ and $\delta_1$ such that 
\[\sup_{s \in G} \frac{1}{m(K_{\delta})} \int_{K_{\delta}} \psi(x+s)dx < \alpha + \varepsilon \]
whenever $\delta \ge \delta_0$ and 
\[\inf_{s \in G} \frac{1}{m(K_{\delta})} \int_{K_{\delta}} \psi(x+s)dx > \alpha - \varepsilon \]
whenever $\delta \ge \delta_1$. Hence, let $\delta_2 = \max(\delta_0, \delta_1)$ and we have
\[-\varepsilon < \frac{1}{m(K_{\delta})} \int_{K_{\delta}} \psi(x+s)dx - \alpha < \varepsilon \]
for every $s \in G$ whenever $\delta \ge \delta_2$. This means uniform convergence of the net in question. This completes the proof.
\end{prf}

Now, we extend the above result to complex valued functions. 
\begin{thm}
Let $\psi \in L^{\infty}(G)$. Then, $\psi$ is alomst convergent to $\alpha$ if and only if 
\[\lim_{\delta} \frac{1}{m(K_{\delta})} \int_{K_{\delta}} \psi(x+s)dx = \alpha \]
uniformly in $s \in G$. 
\end{thm}

\begin{prf}
Let us denote $\psi(x) = u(x) + iv(x)$, where $u, v$ in $L^{\infty}_{\mathbb{R}}(G)$ are real and imaginary parts of $\psi$ respectively and let $\alpha = \beta + i\tau$ be a complex number. Then, 
\[\varphi(\psi) = \alpha = \beta + i\tau \]
for every $\varphi \in \mathcal{T}(G)$ if and only if 
\[\varphi(u) = \beta, \quad \varphi(v) = \tau \]
for every $\varphi \in \mathcal{T}(G)$. By Theorem $2.3$, this is equivalent to 
\[\lim_{\delta} \frac{1}{m(K_{\delta})} \int_{K_{\delta}} u(x+s)dx = \beta, \quad  \lim_{\delta} \frac{1}{m(K_{\delta})} \int_{K_{\delta}} v(x+s)dx = \gamma \]
uniformly in $s \in G$, which is equivalent to the assertion of the theorem. 
\end{prf}

We give examples of summing nets. It can be easily varified that 
\[K_k = [-k, k] \ (k \in \mathbb{N}), \]
\[K_{\theta} = [-\theta, \theta] \ (\theta > 0) \]
are summing nets for $G = \mathbb{Z}$ and $\mathbb{R}$, respectively. Then, we have the following formulation of almost convergence on $\mathbb{Z}$ and $\mathbb{R}$:

$\psi \in L^{\infty}(\mathbb{Z})$ is almost convergent to $\alpha$ if and only if 
\[\lim_{k \to \infty} \frac{1}{2k+1} \sum_{i =n-k}^{n+k} \psi(i) = \alpha \]
uniformly in $n \in \mathbb{Z}$. Similarly, $\psi \in L^{\infty}(\mathbb{R})$ is almost convergent to $\alpha$ if and only if 
\[\lim_{\theta \to \infty} \frac{1}{2\theta} \int_{x-\theta}^{x+\theta} \psi(t)dt = \alpha \]
uniformly in $x \in \mathbb{R}$.

We note that the former is a two-sided version of Sucheston's result (Theorem 1.1) and the latter was given by Raimi \cite{Rai1} (see also \cite{Rai2}). Another kind of analytic expression for almost convergence on $\mathbb{R}$ was given in \cite{Kuni}.

\section{Tauberian theorems concerning almost convergence}
In this section, we introduce two other kinds of convergence on locally compact abelian groups and establish the relationship between those convergence and almost convergence. 

Let $\psi \in L^{\infty}(G)$. We write $\psi \xrightarrow{w^*c} \alpha$ if and only if $w^*\mathchar`-\lim_{x \to \infty} \psi_x = \alpha$ where the symbol $w^*\mathchar`-\lim$ denotes the weak* convergence in $L^{\infty}(G)$. Namely, $\psi \xrightarrow{w^*c} \alpha$ if and only if 
\[\lim_{x \to \infty} \int f(t)\psi(x-t)dt = \alpha \int_G f(t)dt \]
holds for every $f \in L^1(G)$. We write $\psi \xrightarrow{c} \alpha$ if and only if for any $\varepsilon > 0$ there exists a compact subset $C$ of $G$ such that $|\psi - \alpha| < \varepsilon$ whenever $x \notin C$. Note that by the celebrated Wiener's Tauberian theorem, for $\psi \xrightarrow{w^*c} \alpha$ to be hold, it is sufficient that the above equation holds for a single $f \in L^1(G)$ such that $\hat{f}$ does not vanish on $\Gamma$ (see \cite{Rud}).

We can easily see that the following implication holds:
\[\psi \xrightarrow{c} \alpha \Rightarrow \psi \xrightarrow{w^*c} \alpha \Rightarrow \psi \xrightarrow{ac} \alpha \]
In fact, the first half is obvious. Assume that $\psi \xrightarrow{w^*c} \alpha$. This means that for any fixed $f \in P(G)$, we have $f * \psi \xrightarrow{c} \alpha$. By Lemma 2.1 and the proof of Theorem 2.3, $\psi \xrightarrow{ac} \alpha$ if and only if $f * \psi \xrightarrow{ac} \alpha$. Since $f*\psi \xrightarrow{c} \alpha$ obviously implies $f*\psi \xrightarrow{ac} \alpha$, we obtain $\psi \xrightarrow{ac} \alpha$.

We now consider Tauberian theorems. Namely, conditional inverses of the above implication. We say that $\psi \in L^{\infty}(G)$ is slowly oscillating if for any $\varepsilon > 0$, there exist a compact neighborhood $U$ of $0$ and a compact subset $C$ of $G$ such that $|\psi(x) - \psi(y)| < \varepsilon$ whenever $x-y \in U$ and $x \notin C$. Then we have the following result (see \cite{Rud} for proof):
\begin{thm}
Let $\psi \in L^{\infty}(G)$. If $\psi \xrightarrow{w^*c} \alpha$ and $\psi$ is slowly oscillating, then $\psi \xrightarrow{c} \alpha$.
\end{thm}

Now we give a Tauberian theorem which induces weak* convergence from almost convergence.
\begin{thm}
Let $\psi \in L^{\infty}(G)$. If $\psi \xrightarrow{ac} \alpha$ and $\psi-\psi_s \xrightarrow{w^*c} 0$ for every $s \in G$, then $\psi \xrightarrow{w^*c} \alpha$.
\end{thm}

\begin{prf}
First, we show the assertion for an element of $C_{bu}(G)$, in which case the condition $\psi-\psi_s \xrightarrow{w^*c} 0$ is equivalent to $\psi-\psi_s \xrightarrow{c} 0$. Assume that $\psi \in C_{bu}(G)$ and $\psi \xrightarrow{ac} \alpha$ and $\psi-\psi_s \xrightarrow{c} 0$ for every $s \in G$. Fix a positive number $\varepsilon > 0$. Let $K_{\delta_0}$ be an element of a summing net $\{K_{\delta}\}$ for $G$ such that 
\[\left|\frac{1}{m(K_{\delta_0})} \int_{K_{\delta_0}} \psi_x(t)dt - \alpha\right| < \frac{\varepsilon}{2} \]
holds for every $x \in G$. Let $U$ be an open neighborhood of $0$ such that $x-y \in U$ implies $|\psi(x)-\psi(y)| < \frac{\varepsilon}{4}$. We consider the family of opensubsets $\{U_x\}_{x \in K_{\delta_0}}$, where $U_x = x+U = \{x+y : y \in U\}$. This is an open covering of $K_{\delta_0}$ and thus, we can select finite open subcovering $\{U_{x_i}\}_{i=1}^n$ of $K_{\delta_0}$ since $K_{\delta_0}$ is compact:
\[K_{\delta_0} \subseteq \bigcup_{i=1}^n U_{x_i}. \]
For each $x_i \; (1 \le i \le n)$, choose a compact subset $C_i$ of $G$ such that 
\[x \notin C_i \Rightarrow |\psi(x) - \psi(x+x_i)| < \frac{\varepsilon}{4}. \]
Let $C = \cup_{i=1}^n C_i$. Then, we have 
\[x \notin C \Rightarrow |\psi(x) - \psi(x+y)| < \frac{\varepsilon}{2} \]
for every $y \in K_{\delta_0}$. In fact, assume that $y$ is in $U_{x_j}$, then we have
\begin{align}
|\psi(x) - \psi(x+y)| &= |\psi(x) - \psi(x+x_j) + \psi(x+x_j) - \psi(x+y)| \notag \\
&\le |\psi(x) - \psi(x+x_j)| + |\psi(x+x_j) - \psi(x+y)| \notag  \\
&< \frac{\varepsilon}{4} + \frac{\varepsilon}{4} = \frac{\varepsilon}{2}. \notag
\end{align}
Now observe that for any $x \notin C$, we have
\begin{align}
|\psi(x) - \alpha| &\le \left|\frac{1}{m(K_{\delta_0})} \int_{K_{\delta_0}} \psi(x+t)dt - \psi(x)\right| + \left|\frac{1}{m(K_{\delta_0})} \int_{K_{\delta_0}} \psi(x+t)dt-\alpha\right| \notag \\
&\le \left|\frac{1}{m(K_{\delta_0})} \int_{K_{\delta_0}} \psi(x+t)dt - \frac{1}{m(K_{\delta_0})} \int_{K_{\delta_0}} \psi(x)dt\right| + \frac{\varepsilon}{2} \notag \\
&\le \frac{1}{m(K_{\delta_0})} \int_{K_{\delta_0}} |\psi(x+t)-\psi(x)|dt + \frac{\varepsilon}{2} \notag \\
&= \varepsilon. \notag 
\end{align}
This means that $\psi \xrightarrow{c} \alpha$. 

Next we consider a general case. Let $\psi \in L^{\infty}(G)$. Take a function $f \in P(G)$ such that $\hat{f}$ does not vanish on $\Gamma$ and consider $f * \psi \in C_{bu}(G)$. Then, by $\psi \xrightarrow{ac} \alpha$, we have $f*\psi \xrightarrow{ac} \alpha$. The assumption $\psi - \psi_s \xrightarrow{w^*c} 0$ implies that $(f * \psi) - (f*\psi)_s \xrightarrow{c} 0$. Hence, by the above result, we obtain $f * \psi \xrightarrow{c} \alpha$ which in turn implies $\psi \xrightarrow{w^*c} \alpha$ by Wiener's Tauberian theorem. We complete the proof.
\end{prf}

\begin{rem}
In Theorems $3.1$ and $3.2$, the conditions that $\psi$ is slowly oscillating and $\psi-\psi_s \xrightarrow{w^*c} 0$ for every $s \in G$ are also sufficient condition to $\psi \xrightarrow{c} \alpha$ and $\psi \xrightarrow{w^*c} \alpha$, respectively.
\end{rem}
 
\section{Spectral synthesis of bounded measurable functions}
Let $G$ be a locallly compact abelian group and $\Gamma$ be its dual group. An element $\lambda$ in $\Gamma$ is denoted as $\chi_{\lambda}(x)$ when it is viewed as a character on $G$. For a function $f$ in $L^1(G)$, its Fourier transform $\hat{f}$ is defined by
\[\hat{f}(\lambda) = \int_G f(x)\chi_{\lambda}(-x)dm(x), \quad \lambda \in \Gamma. \]
The purpose of this section is to develop the spectrum theory of $C_{bu}(G)$ and $C_{bu}(G)^*$. For the following contents, we refer the readers to \cite{Rud} as a basic literature.

For $f \in L^1(G)$, let $Z(f)$ be the zero set of the Fourier transform of $f$. 
Let $I$ be a closed ideal of $L^1(G)$. We denote by $Z(I)$ the intersection of the zero sets of the Fourier transforms of elements in $I$, namely, 
\[Z(I) = \bigcap_{f \in I} Z(f). \]
By definition, each $Z(I)$ is a closed set of $\Gamma$. Note that, in general, $Z(I)$ does not uniquely determine $I$. In other words, there exist distinct closed ideals $I$ and $I^{\prime}$ such that $Z(I) = Z(I^{\prime})$ holds true. However, there is a closed set $C$ of $\Gamma$ which is the zero set of a unique closed ideal of $L^1(G)$. Such a closed set is called a {\it spectral synthesis set}. The following result is a fundamental result of spectral synthesis thoery.

\begin{thm}
Let $I$ be a closed ideal of $L^1(G)$ and $f \in L^1(G)$. If the internal of $Z(f)$ contains $Z(I)$, then $f \in I$ holds true.
\end{thm}
For a closed set $C$ of $\Gamma$, let $I_1(C)$ be the set of all $f \in L^1(G)$ such that $C \subseteq Z(f)$ and let $I_0(C)$ be the closure of the set of all $f \in L^1(G)$ such that $C \subseteq \text{Int}_{\Gamma} Z(f)$, the interior of $Z(f)$ in $\Gamma$. Then, by definition, $I_1(C)$ is the largest closed ideal of $L^1(G)$ with $Z(I) = C$ and, by Theorem 3.1, $I_0(C)$ is the smallest closed ideal of $L^1(G)$ with $Z(I) = C$. Note that $C$ is a spectral synthesis set if and only if $I_0(C) = I_1(C)$.

Now, following classical theory, we define spectrum of subspaces and elments of $L^{\infty}(G)$. Let $\Phi$ be a weak* closed invariant subspace of $L^{\infty}(G)$, namely, $\Phi$ is a subspace closed with respect to weak* topology of $L^{\infty}(G)$ and $\psi \in \Phi$ implies $\psi_s \in \Phi$ for every $s \in G$. Related to $\Phi$, we define the closed ideal $J(\Phi)$ of $L^1(G)$ by the set of all functions $f$ in $L^1(G)$ such that $f * \psi = 0$ for every $\psi \in \Phi$:
\[J(\Phi) = \{f \in L^1(G) : f * \psi = 0 \ \forall \psi \in \Phi\}. \]
Then, we define the spectrum $\mathrm{sp}(\Phi)$ of $\Phi$ by $Z(J(\Phi))$.

In the same way, we define spectrum of each member of $L^{\infty}(G)$. For each $\psi$ in $L^{\infty}(G)$, let $J(\psi)$ be the set of all functions $f$ in $L^1(G)$ such that $f * \psi = 0$. Then, we define the spectrum $\mathrm{sp}(\psi)$ of $\psi$ by $Z(J(\psi))$. It is easy to confirm that $\mathrm{sp}(\psi) = \mathrm{sp}(\Phi(\psi))$, where $\Phi(\psi)$ is the weak* closed invariant subspace of $L^{\infty}(G)$ generated by $\psi$.

We give another description of spectrum of $\Phi$. Let $\langle X, X^* \rangle$ be any dual pair of locally convex spaces. For subspaces $E$ of $X$ and $E^*$ of $X^*$, let us define their annihilators as follows:
\[E^{\bot} = \{\varphi \in X^* : \varphi(x) = 0 \ \forall x \in E\}, \]
\[E^{*\bot} = \{x \in X : \varphi(x) = 0 \ \forall \varphi \in E^* \}. \]
Note that $E^{\bot}$ is a weak* closed subspace of $X^*$ and $E^{*\bot}$ is a closed subspace of $X$. By the Hahn-Banach theorem, 
\[(E^{\bot})^{\bot} = E, \quad ({E^*}^{\bot})^{\bot} = E^* \]
holds whenever $E$ and $E^*$ are closed subspaces of $X$ and $X^*$ respectively. The important fact about the dual pair $\langle L^1(G), L^{\infty}(G) \rangle$ is the correspondence between the closed ideals of $L^1(G)$ and weak* invariant subspaces of $L^{\infty}(G)$; let $I$ be a closed ideal of $L^1(G)$, then its annihilator $\Phi = I^{\bot}$ is a closed invariant subspace of $L^{\infty}(G)$ and vice versa. 

\begin{thm}
Let $\Phi$ be a weak* closed invariant subspace of $L^{\infty}(G)$. Then, $\mathrm{sp}(\Phi)$ is the set of all chracters contained in $\Phi$ :
\[\mathrm{sp}(\Phi) = \{\lambda \in \Gamma : \chi_{\lambda} \in \Phi\}. \]
\end{thm}

\begin{prf}
First, suppose $\chi_{\lambda} \in \Phi$. Let $f$ be in $J(\Phi)$. Then, in particular, we have
\begin{align}
f * \chi_{\lambda}(x) &= \int_G f(t)\chi_{\lambda}(x-t)dm(t) \notag \\
&= \int_G f(t)\chi_{\lambda}(x)\chi_{\lambda}(-t)dm(t) \notag \\
&= \chi_{\lambda}(x)\hat{f}(\lambda) = 0 \notag 
\end{align}
for every $x \in G$. Thus, $\hat{f}(\lambda) = 0$ and $\lambda \in Z(J(\Phi)) = \mathrm{sp}(\Phi)$.

On the other hand, let us assume that $\chi_{\lambda} \not\in \Phi$. Hence, by the Hahn-Banach theorem, there exists some $f \in L^1(G)$ such that $\langle f^*, \psi \rangle = 0$ for every $\psi \in \Phi$ and $\langle f^*, \chi_{\lambda} \rangle = 1$, where $f^*(x) = f(-x)$. That is, 
\[\int_G f^*(t)\psi(t)dt = \int_G f(-t)\psi(t)dm(t) = f * \psi(0) = 0 \]
for each $\psi \in \Phi$ and 
\[\int_G f^*(t)\chi_{\lambda}(t)dt = \int_G f(-t)\chi_{\lambda}(t)dm(t) = \hat{f}(\lambda) = 1 \]
holds true. Since $\Phi$ is translation invariant, the first equation also holds for any translate $\psi_x$ of $\psi$ and we obtain $f * \psi = 0$. Hence, $f$ is in $J(\Phi)$ and $\hat{f}(\lambda) = 1$, it holds that $\lambda \not\in \mathrm{sp}(\Phi)$. This completes the proof.
\end{prf}

 Obviously, the smallest weak* closed subspace of $L^{\infty}(G)$ with $\mathrm{sp}(\Phi) = C$ is the subspace $\Phi_1(C)$ which is generated by the characters $\{\chi_{\lambda}\}_{\lambda \in C}$. Note that if $I = \Phi^{\bot}$, or, equivalently, $I^{\bot} = \Phi$, $\mathrm{sp}(\Phi) = Z(I)$ holds true. This implies that $\Phi_1(C) = I_1(C)^{\bot}$. Further, let $\Phi_0(C)$ be the annihilaor of $I_0(C)$, that is, $\Phi_0(C) = I_0(C)^{\bot}$, then $\Phi_0(C)$ is the largest weak* closed invariant subspace with $\mathrm{sp}(\Phi) = C$. We have $\Phi_0(C) = \Phi_1(C)$ if and only if $\mathrm{sp}(\Phi)$ is a spectral synthesis set. This observation implies the following assertion:
\begin{thm}
Let $\Phi$ be a weak* closed invariant subspace of $L^{\infty}(G)$. If $\mathrm{sp}(\Phi)$ is a spectral synthesis set, then $\Phi$ is synthesized by its spectrum. Namely, each member $\psi$ of $\Phi$ is a weak* limit of a net of trigonometric polynomials
\[\psi_{\alpha}(x) = \sum_{k=1}^{n_{\alpha}} c_k^{\alpha} e^{it_k^{\alpha}x} \]
where $c_k^{\alpha} \in \mathbb{C}$ and $t_k^{\alpha} \in \mathrm{sp}(\Phi)$ for every $\alpha$ and $1 \le k \le n_{\alpha}$.
\end{thm}

Following \cite{Fore}, we can also define spectrum for elements of $C_{bu}(G)^*$. Let $\varphi \in C_{bu}(G)^*$ and $f \in L^1(G)$. We define their convolution $f * \varphi$, which is also an element of $C_{bu}(G)^*$, by
\[f * \varphi(\psi) = \varphi(f * \psi), \quad \psi \in L^{\infty}(G). \]
For each $\varphi \in C_{bu}(G)^*$, let $J(\varphi)$ be the closed ideal of $L^1(G)$ consisting of those elements $f$ for which $f * \varphi = 0$. Then, we define the spectrum sp$(\varphi)$ of $\varphi$ by $Z(J(\varphi))$.

The following lemma is basic for arguments what follows.
\begin{lem}
Let $\psi$ be in $L^{\infty}(G)$, $\varphi$ be in $C_u(\mathbb{R})^*$ and $f$ be in $L^1(\mathbb{R})$. Then, the following results holds true.  \\
$\mathrm{(i)}$ $\mathrm{sp}(f*\psi) \subseteq \mathrm{sp}(\psi) \cap \mathrm{supp} \; \hat{f}$, \\
$\mathrm{(ii)}$ $\mathrm{sp}(f*\varphi) \subseteq \mathrm{sp}(\varphi) \cap \mathrm{supp} \; \hat{f}$, \\
where $\mathrm{supp} \; \hat{f}$ denotes the support of $\hat{f}$.
\end{lem}

\begin{prf}
$\mathrm{(i)}$ First, note that the inclusion $\mathrm{sp}(f * \psi) \subseteq \mathrm{sp}(\psi)$ is obvious by definition of spectrum. Fix arbitrary $\lambda \not\in \mathrm{supp} \; \hat{f}$. Let $g$ be a funciton in $L^1(G)$ such that $\hat{g} = 0$ on a neighborhood of $\mathrm{supp} \; \hat{f}$ and $\hat{g}(\lambda) = 1$. Then, note that ${(f * g)}^{\widehat{}}(\lambda) = \hat{f}(\lambda)\hat{g}(\lambda) = 0$ for every $\lambda \in \Gamma$, which means that $f * g = 0$ by the uniqueness theorem of the Fourier transform. Hence, we obtain $g * (f * \psi) = (g * f) * \psi = 0$ and thus, $g$ is in $J(f*\psi)$. Since $\lambda \not\in Z(g)$, we conclude that $\lambda \notin \mathrm{sp}(f * \psi)$. \\
$\mathrm{(ii)}$ This assertion can be proved in the same way as $\mathrm{(i)}$ and we omit the proof.
\end{prf}

\begin{lem}
Let $\psi$ be in $L^{\infty}(G)$ and $\varphi$ be in $C_{bu}(G)^*$. Then, we have the following results. \\
$\mathrm{(i)}$ If $\mathrm{sp}(\psi) = \emptyset$, then $\psi = 0$, \\
$\mathrm{(ii)}$ If $\mathrm{sp}(\varphi) = \emptyset$, then $\varphi = 0$.
\end{lem}

\begin{prf}
$\mathrm{(i)}$ Suppose that $\mathrm{sp}(\psi) = \emptyset$ holds. By definition of spectrum and Wiener's Tauberian theorem, which asserts that if a closed ideal $I$ in $L^1(G)$ satisfies $Z(I) = \emptyset$, then $I = L^1(G)$, we conclude that $f * \psi = 0$ for every $f$ in $L^1(G)$. This means $\psi = 0$ by the Hahn-Banach theorem.

\noindent
$\mathrm{(ii)}$ Suppose that $\mathrm{sp}(\varphi) = \emptyset$, that is, $f * \varphi(\psi) = \varphi(f*\psi) = 0$ for every $\psi \in C_u(G)$. Take an approximate identity $\{f_{\alpha}\}$ of $L^1(G)$, observe that
\[\lim_{\alpha} \|f_{\alpha} * \psi - \psi\|_{\infty} = 0 \]
holds for every $\psi \in C_u(G)$. Then, for any $\psi \in C_u(G)$, we have
\begin{align}
|\varphi(\psi)| &\le |\varphi(\psi - f_{\alpha} * \psi)| + |\varphi(f_{\alpha} * \psi)| \notag \\
&= |\varphi(\psi - f_{\alpha} * \psi)|, \notag 
\end{align}
which tends to $0$ as $\alpha \to \infty$. Thus, it follows that $\varphi(\psi) = 0$, completing the proof.  
\end{prf}

We define the subspaces of $C_{bu}(G)$ and $C_{bu}(G)^*$ as follows: 
\[E_A = \{\psi \in C_{bu}(G) : \text{sp}(\psi) \subseteq A\}; \]
\[E^*_A = \{\varphi \in C_{bu}(G)^* : \text{sp}(\varphi) \subseteq A\}, \]
where $A$ is a subset of $\Gamma$. Observe that if $A$ is a closed subset of $\Gamma$, then $E_A = \Phi_0(A) \cap C_{bu}(G)$ holds true and thus, $E_A$ is a closed subspace of $C_{bu}(\mathbb{R})$. For spaces $E^*_A$, we have the following result.

\begin{lem}
Let $A$ be a closed subset of $\Gamma$. Then, $E^*_A$ is a weak*-closed subspace of $C_{bu}(G)^*$. In particular, it is norm closed.
\end{lem}

\begin{prf}
Let $\{\varphi_{\alpha}\}$ be a net of elements of $E^*_A$ such that $w^*\mathchar`-\lim_{\alpha} \varphi_{\alpha} = \varphi$. We show that $\varphi \in E^*_A$, that is, $\mathrm{sp}(\varphi) \subseteq A$. Take arbitrary $\lambda \not\in A$. Then, there exists a function $f$ in $L^1(G)$ such that $\hat{f} = 0$ on a some neighborhood of $A$ and $\hat{f}(\lambda) = 1$. Note that $f * \varphi_{\alpha} = 0$ for every $\alpha$. Then, for any $\psi$ in $L^{\infty}(G)$, we have
\[(f * \varphi)(\psi) = \varphi(f * \psi) = \lim_{\alpha} \varphi_{\alpha} (f * \psi) = \lim_{\alpha} f * \varphi_{\alpha}(\psi)  = 0. \]
Hence, we obtain $f * \varphi = 0$ and thus, $f \in J(\varphi)$. Since $\lambda \not\in Z(f)$, we have $\lambda \not\in \mathrm{sp}(\varphi)$, completing the proof. 
\end{prf}

The following result will be used in Section 8, which states the continuity of spectrum with respect to the weak* topology.
\begin{lem}
Let $A$ be a closed subset of $\Gamma$. Suppose that $\{\psi_{\alpha}\}$ be a net of $L^{\infty}(G)$ with $\mathrm{sp} (\psi_{\alpha}) \subseteq A$ for every $\alpha$. If $w^*\mathchar`-\lim_{\alpha} \psi_{\alpha} = \psi$, namely, $\{\psi_{\alpha}\}$ converge to $\psi$ in the weak* sense, then we have $\mathrm{sp} (\psi) \subseteq A$.
\end{lem}

\begin{prf}
In fact, for each $f \in L^1(G)$ with $A \subseteq {\rm Int}_{\Gamma} Z(f)$, in other words, for $f \in I_0(A)$, we have $f * \psi_{\alpha} = 0$ for every $\alpha$. Thus, we have
\[f * \psi = \lim_{\alpha} f * \psi_{\alpha} = 0, \]
which means that $f \in J(\psi)$. Since $\cap_{f \in I_0(A)} Z(\hat{f}) = A$, we obtain $\mathrm{sp}(\psi) = \cap_{f \in J(\psi)} Z(\hat{f}) \subseteq A$. This complets the proof.
\end{prf}

In what follows, we give a description of the spaces $E^*_A$ in terms of spaces $E_A$. Note that we consider the dual pair $\langle C_{bu}(G), C_{bu}(G)^* \rangle$.
\begin{thm}
Let $A$ be a closed subset of $\Gamma$. Then, we have
\[{E^*_A}^{\bot} = {\rm cl}_{C_{bu}(G)} E_{A^c}, \]
where the symbol ${\rm cl}_{C_{bu}(G)} E$ denotes the closure of a subspace $E$ in $C_{bu}(G)$ $($with respect to the norm topology$)$.
\end{thm}

\begin{prf}
For each $\psi \in C_{bu}(G)$ and $\varphi \in C_{bu}(G)^*$, Let us define $\psi^{\prime}(x) \in C_{bu}(G)$ by $\psi^{\prime}(x) = \varphi(\psi_x)$. For any $f \in L^1(G)$, by the fact mentioned in the proof of Theorem $2.1$, we have
\[
\varphi(f * \psi) = \int_G \varphi(\psi_{-t})f(t)dt. 
\]
Replacing $\psi$ by $\psi_x$ where $x \in G$, we obtain 
\begin{equation}
f * \varphi(\psi_x) = \varphi(f * \psi_x) = \int_G \varphi(\psi_{x-t})f(t)dt = \int_G \psi^{\prime}(x-t)f(t)dt = f * \psi^{\prime}(x). 
\end{equation}
From this equation, we can deduce the following result:
\begin{equation}
\mathrm{sp}(\psi^{\prime}) \subseteq \mathrm{sp}(\psi) \cap \mathrm{sp}(\varphi). 
\end{equation}
In fact, assume that $f \in J(\psi)$, that is, $f * \psi = 0$. Then, for any $x \in G$, it holds that $f * \psi_x = (f * \psi)_x = 0$. Thus, by the equation $(3)$, if $f \in J(\psi)$, we have $f * \psi^{\prime} = 0$, that is, $f \in J(\psi^{\prime})$. Hence, we obtain $\mathrm{sp}(\psi^{\prime}) \subseteq \mathrm{sp}(\psi)$ since $J(\psi) \subseteq J(\psi^{\prime})$.

Next, assume that $f \in J(\varphi)$, that is, $f * \varphi = 0$. Again, by equation $(3)$, we obtain $f * \psi^{\prime} = 0$, that is, $f \in J(\psi^{\prime})$. This shows that $\mathrm{sp}(\psi^{\prime}) \subseteq \mathrm{sp}(\varphi)$. We have obtained the desired relation $(4)$.

Now suppose that $\varphi \in E_A^*(G)$ and $\psi \in E_{A^c}(G)$. Then, by $(4)$, we have $\mathrm{sp}(\psi^{\prime}) \subseteq \mathrm{sp}(\psi) \cap \mathrm{sp}(\varphi) \subseteq A \cap A^c = \emptyset$. Thus, $\mathrm{sp}(\psi^{\prime}) = \emptyset$ holds true. Then, by Lemma $4.2 \; (1)$, we have $\psi^{\prime} = 0$ and $\varphi(\psi) = 0$, that is, $\psi \in {E^*_A}^{\bot}$. Now we obtain
\[E_{A^c} \subseteq {E^*_A}^{\bot}, \]
which means that
\[{\rm cl}_{C_{bu}(G)} E_{A^c}  \subseteq {E^*_A}^{\bot}. \]
We now show the reverse inclusion. Let us assume that $\varphi(\psi) = 0$ for every $\psi \in E_{A^c}$. We show that $\mathrm{sp}(\varphi) \subseteq A$. Let $\lambda \not\in A$ and take $f \in L^1(G)$ such that $\hat{f} = 0$ on some neighborhood $U$ of $A$ with $\lambda \not\in U$ and $\hat{f}(\lambda) = 1$. Then, $f * \varphi(\psi) = \varphi(f * \psi) = 0$ for every $\psi \in C_{bu}(G)$ since $\mathrm{sp}(f * \psi) \subseteq \mathrm{supp} \hat{f} \subseteq U^c \subseteq A^c$. Thus, $f$ is a function in $J(\varphi)$ with $\hat{f}(\lambda) = 1$ and we conclude that $\lambda \not\in \mathrm{sp}(\varphi)$. Thus, we have shown that
\[{E_{A^c}}^{\bot} \subseteq E^*_A. \]
Considering the annihilators of the both sides of the above equation and then taking the closure, we obtain the inclusion relation
\[{E^*_A}^{\bot} \subseteq {\rm cl}_{C_{bu}(G)} E_{A^c}. \]
This completes the proof.

\end{prf}

\section{An application to topologically invariant means}
We apply Theorem 4.4 to the case $C = \{0\}$, that is, the singleton consisting of only the unit of $G$. 
The following result is essentially due to Muhly \cite{Muh}, who proved the special case of $G = \mathbb{R}$.

\begin{thm}
For $\varphi \in C_{bu}(G)^*$, $\varphi$ is an invariant mean if and only if $\text{sp}(\varphi) = \{0\}$.
\end{thm}

\begin{prf}
First, assume that $\varphi$ is invariant. Then, for any $f \in L^1(G)$ for which $Z(f)$ contains $\{0\}$, that is, $\int_G f(x)dm(x) = 0$, we have
\begin{align}
f * \varphi(\psi) &= \varphi(f * \psi) = \int_G \varphi(\psi_{-t})f(t)dt \notag \\
&= \int_G \varphi(\psi)f(t)dt = \varphi(\psi) \int_G f(t)dt = 0. \notag
\end{align}
Therefore, we conclude that $\mathrm{sp}(\varphi) \subseteq \{0\}$. Note that if $\mathrm{sp}(\varphi) = \emptyset $, by Lemma $4.2$, $\varphi = 0$ holds true. Hence, we conclude that $\mathrm{sp}(\varphi) = \{0\}$.

Conversely, assume that $\text{sp}(\varphi) = \{0\}$. Then, for any $f \in I_0(\{0\})$ and $\psi \in C_{bu}(G)$, we have
\[f * \varphi(\psi) = \int_{\mathbb{R}} \varphi(\psi_{-t})f(t)dt = 0. \]
For any $x \in G$, set $\psi^{\prime}(x) = \varphi(\psi_x)$ and we have
\[f * \varphi(\psi_x) = \int_{\mathbb{R}} \varphi(\psi_{x-t})f(t)dt = \int_{\mathbb{R}} \psi^{\prime}(x-t)f(t)dt = 0. \]
This means that $\text{sp}(\psi^{\prime}) = \{0\}$. Since $\{0\}$ is a spectral synthesis set (see \cite{Rud}), we have $\Phi(\psi^{\prime}) = \Phi_1(\{0\})$ = $\mathbb{R}$. Hence, $\psi^{\prime}(x)$ is a constant function, which means that $\varphi$ is translation invariant. This completes the proof. 
\end{prf}

Let $T(G)$ be the closed subspace of $L^{\infty}(G)^*$ spaned by $\mathcal{T}(G)$. We determine the annihilator of $T(G)$ via Theorem 4.4. For any subset $A$ of $\Gamma$, let us define the subspace $E^{\prime}_A$ of $L^{\infty}(G)$ by
\[E^{\prime}_A = \{\psi \in L^{\infty}(G) : \text{sp}(\psi) \subseteq A\}.  \]
\begin{thm}
For a locally compact abelian group $G$, the following assertion holds true.
\[T(G)^{\bot} = cl_{L^{\infty}(G)} E^{\prime}_{G \setminus \{0\}}. \]
\end{thm}

\begin{prf}
Let $\varphi \in T(G)$. Fix $\psi \in E^{\prime}_{G \setminus \{0\}}$. For any $f \in P(G)$, we have $f * \psi \in C_{bu}(G)$ and $\text{sp}(f * \psi) \subseteq G \setminus \{0\}$ by Lemma $4.1$. Hence, for each $\varphi \in T(G)$, we have
\[\varphi(\psi) = \varphi(f * \psi) = 0 \]
by Theorems $4.4$ and $5.1$. Now, we obtain that
\[T(G)^{\bot} \supseteq cl_{L^{\infty}(G)} E^{\prime}_{G \setminus \{0\}}. \]
We show the reverse inclusion. Suppose $\varphi \in L^{\infty}(G)$ vanishes on $E^{\prime}_{G \setminus \{0\}}$. We show that $\varphi$ is in $T(G)$, that is, $\varphi(\psi - f * \psi) = 0$ for every $f \in P(G)$ and $\psi \in L^{\infty}(G)$. First, note that for a function $g$ in $L^1(G)$ such that $\hat{g} = 1$ in some neighborhood $U$ of $0$, we have $\varphi(\psi - g * \psi) = 0$ for every $\psi$ in $L^{\infty}(G)$. In fact, to show this, it is sufficient to show that $\psi - g * \psi$ is in $E^{\prime}_{G \setminus \{0\}}$, in other  words, $\mathrm{sp}(\psi - g * \psi) \subseteq G \setminus \{0\}$. Take a function $h$ in $L^1(G)$ such that $\hat{h} = 0$ outside $U$ and $\hat{h}(0) = 1$. Observe that
\[h * (\psi - g * \psi) = (h - g * h) * \psi. \]
Here, we have $(h-g*h)^{\widehat{}}(\lambda) = \hat{h}(\lambda) (1 - \hat{g}(\lambda)) = 0$ for every $\lambda \in \Gamma$. Hence, $h - g * h = 0$ and thus $h \in J(\psi-g*\psi)$. Since $\hat{h}(0) = 1$, we conclude that $0 \not\in \mathrm{sp}(\psi - g*\psi)$. We obtain the desired assertion. 

For any $f \in P(G)$ and $\varepsilon > 0$, there exists a function $h$ in $L^1(G)$ such that  $\hat{h} = 1$ on some neighborhood of $0$ and $\|f - h\|_1 < \varepsilon$(See \cite{Rud}, Theorem $2.6.5$). Hence, we have
\begin{align}
|\varphi(\psi - f * \psi)| &= |\varphi(\psi - h * \psi) + \varphi(h * \psi - f * \psi)| \notag \\
&\le |\varphi(\psi-h*\psi)| + |\varphi(h*\psi-f*\psi)| \notag \\
&= |\varphi(h*\psi-f*\psi)| \notag \\ 
&= |\varphi((h-f) * \psi))| \notag \\
&\le \|\varphi\| \|h-f\|_1 \|\psi\|_{\infty} \notag \\
&\le \|\varphi\| \|\psi\|_{\infty} \varepsilon. \notag
\end{align}
Since $\varepsilon > 0$ can be arbitrary, we obtain $\varphi(\psi - f *\psi) = 0$, which means that $\psi$ is in $T(G)$. This completes the proof.
\end{prf}

\section{A functional analytic condition for almost convergence}
In this section, applying Theorem $5.2$, we obtain the second necessary and sufficient condition on almost convergence.
Let us define the subspaces of $L^{\infty}(G)$ as follows:
\begin{align}
E(G) &:= \{\psi \in L^{\infty}(G) : \psi \ \text{is almost convergent}\} \notag \\
E_0(G) &:= \{\psi \in L^{\infty}(G) : \psi \ \text{is almost convergent to} \ 0\} \notag
\end{align}
It is clear that the following decomposition holds true:
\[E(G) = \mathbb{R} \oplus E_0(G), \quad \psi \mapsto \alpha + (\psi - \alpha), \]
where $\alpha$ is the number to which $\psi$ almost converges. Theorem 5.2 can be reformulated in terms of almost convergence as follows:

\begin{thm}
Let $G$ be a locally compact abelian group. $\psi \in L^{\infty}(G)$ is almost convergent to $0$ if and only if
\[\psi \in cl_{L^{\infty}(G)} E^{\prime}_{G \setminus \{0\}}. \]
In other words, $\psi$ is almost convergent to $0$ if and only if for any $\varepsilon > 0$, there exists $\psi_1 \in L^{\infty}(G)$ whose spectrum does not intersect some neighborhood of $0$ such that $\|\psi - \psi_1\|_{\infty} < \varepsilon$.
\end{thm}

In what follows, we provide some applications of the theorem. For the proof of the following result, see \cite{Rud}.
\begin{lem}
Let $\mu \in M(\Gamma)$. Let us define $\hat{\mu}^*(x) \in C_{bu}(G)$ by
\[\hat{\mu}^*(x) = \int_G \chi_{\lambda}(x)d\mu(\lambda), \quad x \in G. \]
Then, $\mathrm{sp}(\hat{\mu}^*) = \mathrm{supp} \; \mu$ holds true, where $ \mathrm{supp} \; \mu$ is the support of the measure $\mu$.
\end{lem}

The following result was due to Eberlein \cite{Eber}. He showed that every weakly almost periodic function on $\mathbb{R}$ almost converges and the functions $\hat{\mu}^*$ defined above is weakly almost periodic. Here, we will give a more direct proof using Theorem 6.1.
\begin{thm}
For each $\mu \in M(G)$, $\hat{\mu}^*$ is almost convergent to $\mu(\{0\})$.
\end{thm}

\begin{prf}
Let $\mu_0 = \mu - \mu(\{0\})$. Then, for any $\varepsilon > 0$, there exits a neighborhood $U$ of $0$ such that $|\mu_0|(U) < \varepsilon$. Let us define $\mu_{U^c}(A) = \mu(A \cap U^c)$. Then, $\mathrm{sp}({\hat{\mu}_{U^c}}^*) \cap U = \emptyset$ by Lemma $6.1$ and we have
\[|\hat{\mu}_0^*(x) - {\hat{\mu}_{U^c}}^*(x)| \le \int_U |\chi_{\lambda}(x)|d|\mu|(x) = |\mu_0|(U) < \varepsilon.    \]
Hence, by Theorem $6.1$, we have $\hat{\mu}_0^* \xrightarrow{ac} 0$, which in turn implies that $\hat{\mu}^* \xrightarrow{ac} \mu(\{0\})$, completing the proof.
\end{prf}

\begin{thm}
Let $\{\lambda_n\}_{n=1}^{\infty}$ be a sequence of $\Gamma$. Let $\psi \in L^{\infty}(G)$ be defined by
\[\psi(x) = \sum_{n=1}^{\infty} c_n \chi_{\lambda_n}(x) \ (\text{uniformly bounded and pointwise convergence}). \]
Then, we have $\mathrm{sp}(\psi) \subseteq cl_{\Gamma} {\{\lambda_n\}}_{n=1}^{\infty}$.
\end{thm}

\begin{prf}
Let $\Phi$ be the weak* closed invariant subspace of $L^{\infty}(G)$ generated by $\{\chi_{\lambda_n}\}_{n=1}^{\infty}$. Then, we have obviously $\mathrm{sp}(\Phi) = cl_{\Gamma} \{\lambda_n\}_{n=1}^{\infty}$. Since $\psi \in \Phi$ by the above expression of $\psi$, $\Phi(\psi) \subseteq \Phi$ holds true. Thus, we obtain
\[\mathrm{sp}(\psi) = \mathrm{sp}(\Phi(\psi)) \subseteq \mathrm{sp}(\Phi) = cl_{\Gamma} \{\lambda_n\}_{n=1}^{\infty}. \]
This completes the proof.
\end{prf}

\begin{thm}
Let $\{\lambda_n\}_{n=0}^{\infty} \ (\lambda_0 = 0)$ be a sequence of $\Gamma$ which does not accumulate to $0$. If $\psi \in L^{\infty}(G)$ is expressed by
\[\psi(x) = \sum_{n=0}^{\infty} c_n \chi_{\lambda_n}(x) \ (\text{uniformly bounded and pointwise convergence}). \]
Then, $\psi$ is almost convergent to $c_0$.
\end{thm}

\begin{prf}
By Theorem $6.3$, $\mathrm{sp}(\psi_0) \subseteq cl_{\Gamma} \{\lambda_n\}_{n=1}^{\infty}$, where $\psi_0(x) = \psi(x) - c_0$. By the assumption that $\{\chi_{\lambda_n}\}_{n=1}^{\infty}$ does not accumulate to $0$, there exists some neighborhood $U$ of $0$ such that $U \cap \{\lambda_n\}_{n=1}^{\infty} = \emptyset$. Hence, we have $\mathrm{sp}(\psi_0) \subseteq U^c$ and  by Theorem $6.1$, we have $\psi_0 \xrightarrow{ac} 0$, which shows the desired assertion.

\end{prf}

We give an example in the case $G = \mathbb{R}$. One of the important examples of functions on $\mathbb{R}$ which is expressed by the sum of exponential functions is Dirichlet series.

\begin{cor}
Let $\{a_n\}_{n=1}^{\infty}$ be a sequence of complex numbers and $\psi(s)$ is the Dirichlet series whose coefficients are $\{a_n\}_{n=1}^{\infty}$.
\begin{align}
\psi(s) &= \frac{a_1}{1^s} + \frac{a_2}{2^s} + \cdots \frac{a_n}{n^s} + \cdots \notag \\
&= \frac{a_1}{1^{\sigma}} e^{-it\log1} + \frac{a_2}{2^{\sigma}} e^{-it\log 2} + \cdots \frac{a_n}{n^{\sigma}} e^{-it\log n} + \cdots \ (s = \sigma + it).  \notag 
\end{align}
If $\psi(s)$ converges uniformly boundedly and pointwisely on the line $Re(s)=\sigma$, then, $\psi_{\sigma}(t) = \psi(\sigma + it)$ is almost convergent to $a_1$.
\end{cor}

\section{One-sided almost convergence}
In the following two sections, we consider the special groups of $G = \mathbb{Z}$ and $\mathbb{R}$. For these groups, it is convenient to define almost convergence for functions defined on positive numbers $\mathbb{Z}_+$ and $\mathbb{R}_+ := \{x \in \mathbb{R} : x \ge 0\}$. In what follows, the symbols $G$ stands for $\mathbb{Z}$ or $\mathbb{R}$ and $G_+$ stands for $\mathbb{Z}_+$ or $\mathbb{R}_+$.

Let us define three closed subspaces of $L^{\infty}(G)$ as follows: 
\[L^{\infty}_0(G) := \{\phi \in L^{\infty}(G) : \lim_{|x| \to \infty} \phi(x) = 0\}, \]
\[L^{\infty}_{0,+}(G) := \{\phi \in L^{\infty}(G) : \lim_{x \to \infty} \phi(x) = 0\}, \]
\[L^{\infty}_{0,-}(G) := \{\phi \in L^{\infty}(G) : \lim_{x \to -\infty} \phi(x) = 0\}. \]
For simplicity, we will also use the symbols $L^{\infty}, L^{\infty}_0, L^{\infty}_{0,+}$ and $L^{\infty}_{0,-}$ as abbreviations of $L^{\infty}(G), L^{\infty}_0(G), L^{\infty}_{0,+}(G)$ and $L^{\infty}_{0,-}(G)$, respectively. Now observe that the isomorphism
\begin{align}
L^{\infty}/L^{\infty}_0 &\cong L^{\infty}/(L^{\infty}_{0,+} \cap L^{\infty}_{0,-}) \notag \\
&\cong (L^{\infty}/L^{\infty}_{0,+}) \oplus (L^{\infty}/L^{\infty}_{0,-}) \notag
\end{align}
holds. Here, the image of an equivalence class $[\psi]$ in $L^{\infty}/L^{\infty}_0$ is given by $[\psi_+] + [\psi_-]$, where $\psi_+ := \psi \cdot I_{G_+}$ and $\psi_- := \psi \cdot I_{-G_+}$, where $I_E$ denotes the characteristic function of a subset $E$ of $G$.

Therefore, for any $\varphi$ in $L^{\infty}(G)^*$ with $\varphi = 0$ on $L^{\infty}_0$, we have the decomposition of $\varphi$ given as follows:
\[(L^{\infty}/L^{\infty}_0)^* \cong (L^{\infty}/L^{\infty}_{0,+})^* \oplus (L^{\infty}/L^{\infty}_{0,-})^*, \quad \varphi = \varphi_+ + \varphi_-, \]
where $\varphi_+(\psi) := \varphi(\psi_+)$ and $\varphi_-(\psi) :=\varphi(\psi_-)$ for each $\psi \in L^{\infty}$.
We note that the spaces $(L^{\infty}/L^{\infty}_{0,\pm})^*$ can be identified with the subspace of $L^{\infty}(G)^*$ consisting of those elements vanishing on $L^{\infty}_{0,\pm}$, respectively.

By Theorem $2.2$, for each $\varphi \in \mathcal{T}$, we have $\varphi = 0$ on $L^{\infty}_0$ and hence, $\varphi$ is decomposed into the sum of $\varphi_+ \in (L^{\infty}/L^{\infty}_{0,+})^*$ and $\varphi_- \in (L^{\infty}/L^{\infty}_{0,-})^*$. We use the symbols 
\[\mathcal{T}_+ := \mathcal{T} \cap (L^{\infty}/L^{\infty}_{0,+})^*, \quad \mathcal{T}_- := \mathcal{T} \cap  (L^{\infty}/L^{\infty}_{0,-})^*. \]
Then, $\mathcal{T}$ is equal to $co(\mathcal{T}_+ \cup \mathcal{T}_-)$, the convex hull of $\mathcal{T}_+ \cup \mathcal{T}_-$.

Now we define one-sided almost convergence for elements in $L^{\infty}(G)$. 
\begin{dfn}
We say that $\psi$ is one-sidedly almost convergent to the number $\alpha$ if $\varphi(\psi) = \alpha$ for every $\varphi \in \mathcal{T}_+$.
\end{dfn}
In this case, we write $\psi \overset{\text{oac}}{\rightarrow} \alpha$. We also define one-sided almost convergence for the functions in $L^{\infty}(G_+)$, the space of essentially bounded functions defined on the positive part $G_+$ of $G$, by identifying $\psi$ in $L^{\infty}(G_+)$ as the function on $G$ defined by
\[\tilde{\psi}(x) = 
\begin{cases}
\psi(x),   & x \ge 0, \\
0,          & x < 0.
\end{cases}
\]

\begin{lem}
Let $\psi$ be in $L^{\infty}(G)$ which vanishes on the negative part $-G_+$ of $G$. Then, $\psi \xrightarrow{oac} 0$ if and only if $\psi \xrightarrow{ac} 0$ holds true.
\end{lem}

\begin{prf}
Necessity is obvious. Suppose that $\psi \in L^{\infty}(G)$ is one-sidedly almost convergent to $0$, that is, $\varphi(\psi) = 0$ for every $\varphi \in \mathcal{T}_+$. Since $\psi$ is in $L^{\infty}_{0, -}(G)$, we have $\varphi(\psi) = 0$ for every $\mathcal{T}_-$. Hence, it follows that $\psi \xrightarrow{ac} 0$ by the fact that $\mathcal{T} = co(\mathcal{T}_+ \cup \mathcal{T}_-)$.
\end{prf}

We can provide a similar condition to Theorem $2.3$ for a given $\psi \in L^{\infty}(G)$ to be one-sidedly almost convergent. To this end, we need a following elementary lemma concerning extreme values of means in $(L^{\infty}/L^{\infty}_{0,+})^*$. 
\begin{lem}
Let $\varphi$ be a mean on $L^{\infty}(G)$ which vanishes on $L^{\infty}_{0,+}(G)$. Then, we have
\[\varphi(\psi) \le \sup_{x \in \mathbb{R}_+} \psi(x)  \]
for every $\psi \in L^{\infty}(G)$.
\end{lem}

\begin{thm}
A mean $\varphi$ on $L^{\infty}(G)^*$ is in $\mathcal{T}_+$ if and only if 

\smallskip
\noindent
$\mathrm{(i)} \ G = \mathbb{Z}$ \\
\[\varphi(\psi) \le \overline{p}_+(\psi) := \limsup_{k \to \infty} \sup_{n \in \mathbb{Z}_+} \frac{1}{k} \sum_{i=0}^{k-1} \psi(n+i) \]
holds for every $\psi \in L^{\infty}(\mathbb{Z})$. 

\smallskip
\noindent
$\mathrm{(ii)} \ G = \mathbb{R}$
\[\varphi(\psi) \le \overline{p}_+(\psi) := \limsup_{\theta \to \infty} \sup_{x \in \mathbb{R}_+} \frac{1}{\theta} \int_x^{x+\theta} \psi(t)dt \]
holds for every $\psi \in L^{\infty}(\mathbb{R})$.
\end{thm}

\begin{prf}
Using Lemma $7.2$, necessity can be proved in a similar way to the proof of Theorem $2.1$. For sufficiency, observe that $\overline{p}_+(\psi) \le \overline{p}(\psi)$ for every $\psi \in L^{\infty}(G)$ and thus, by Theorem $2.2$, it follows that a mean $\varphi$ satisfying the condition in the theorem is in $\mathcal{T}$. Furthermore, it is clear that $\overline{p}_+(\psi) = \underline{p}_+(\psi) (:= -\overline{p}_+(-\psi)) = 0$ holds for every $\psi \in L^{\infty}_{0,+}(G)$ and thus, such a $\varphi$ is in $(L^{\infty}/L^{\infty}_{0,+})^*$. As a result, we obtain that $\varphi$ is in $\mathcal{T}_+ = \mathcal{T} \cap (L^{\infty}/L^{\infty}_{0,+})^*$. This completes the proof.
\end{prf}
Now, we can show the following result just as Theorem 2.3 was derived from Theorem 2.2.
\begin{thm}
Let $\psi$ be in $L^{\infty}(G_+)$. Then, $\psi$ is one-sidedly almost convergent to $\alpha$ if and only if 

\smallskip
\noindent
$\mathrm{(i)} \ G = \mathbb{Z}$ 
\[\lim_{k \to \infty} \frac{1}{k} \sum_{i=0}^{k-1} \psi(n+i) = \alpha \]
uniformly in $n \in \mathbb{Z}_+$.

\smallskip
\noindent
$\mathrm{(ii)} \ G = \mathbb{R}$ 
\[\lim_{\theta \to \infty} \frac{1}{\theta} \int_x^{x+\theta} \psi(t)dt = \alpha \]
uniformly in $x \in \mathbb{R}_+$.
\end{thm}
We remark that one-sided almost convergence for $l_{\infty} = L^{\infty}(\mathbb{Z}_+)$ is equivalent to Lorentz's almost convergence and $\mathrm{(i)}$ of the above theorem is the very result Lorentz proved (Theorem $1.1$).

\section{Complex Tauberian Thoerems for almost convergence}
In this section, we study so-called complex Tauberian theorems for almost convergence on the groups $\mathbb{Z}$ and $\mathbb{R}$. One of the main theorems of this section (Theorem $8.6$) can be viewed as an anlogue of the celebrated Wiener-Ikehara theorem. First, following \cite{Kor}, we introduce the following definition.

\begin{dfn}
We say that a function $f(z)$ defined on the unit disc $\mathbb{D} := \{z \in \mathbb{C} : |z| < 1\}$ has $H^1$ boundary behavior at the point $z_0 = e^{it_0} \ (t_0 \in [0, 2\pi])$ if there exist a number $\delta > 0$ and a function $F(e^{it})$ in $L^1(t_0-\delta, t_0 + \delta)$ such that $f(re^{it})$ converges to $F(e^{it})$ in $L^1(t_0-\delta, t_0 + \delta)$ as $r \to 1^-$.
\end{dfn}

\begin{thm}
Let $\psi$ be in $L^{\infty}(\mathbb{Z})$ which vanishes on the negative integers and denote $a_n = \psi(n)$ for $n \ge 0$. Let $\hat{\psi}(z)$ be the function on the unit disc $\mathbb{D}$ defined by
\[\hat{\psi}(z) = \sum_{n=0}^{\infty} a_nz^n. \]
If $\hat{\psi}$ has $H^1$ boundary behaviour at $z = 1$, then $\psi$ is almost convergent to $0$.
\end{thm}

\begin{prf}
For each numbers $r$ with $0 < r < 1$, define the functions $\psi_r(n)$ in $L^1(\mathbb{Z})$ and $\hat{\psi}_r(\theta)$ in $L^1(\mathbb{T})$ by
\[\psi_r(n) := 
\begin{cases}
a_nr^n  & (n \ge 0), \\
0 & (n < 0).
\end{cases}
\]
\[\hat{\psi}_r(\theta) := \hat{\psi}(re^{i\theta}) = \sum_{n=0}^{\infty} a_nr^n e^{in\theta} \ (\theta \in \mathbb{T}), \]
respectively. Let $\delta_0 > 0$ be a number such that
\[\lim_{r \to 1^-} \int_{-\delta_0}^{\delta_0} |\hat{\psi}_r(\theta) - \hat{\psi}(\theta)|\frac{d\theta}{2\pi} = 0 \]
holds true for some $\hat{\psi} \in L^1(-\delta_0, \delta_0)$. Fix $\varepsilon > 0$ and choose a number $0 < \delta < \delta_0$ such that 
\[\int_{-\delta}^{\delta} |\hat{\psi}_r(\theta)| \frac{d\theta}{2\pi} < \varepsilon \]
for every $0 < r <1$. Let $C_{\delta}$ be a function in $L^1(\mathbb{Z})$ such that $\hat{C}_{\delta} = 1$ on $(-\frac{\delta}{2}, \frac{\delta}{2}), \ \hat{C}_{\delta} = 0$ outside $[-\delta, \delta]$ and $0 \le \hat{C}_{\delta} \le 1$ on $\mathbb{T}$. For each $0 < r < 1$, put
\[\hat{\psi}_{r, \delta}(\theta) = \hat{\psi}_r(\theta)(1-\hat{C}_{\delta}(\theta)), \]
and
\[\psi_{r, \delta}(n) = (\hat{\psi}_{r, \delta})^{\widehat{}}(n), \quad n \in \mathbb{Z}. \]
Then, we have
\[\text{sp}(\psi_{r, \delta}) \cap \left(-\frac{\delta}{2}, \frac{\delta}{2}\right) = \emptyset \ (0 < r <1). \]
In fact, by Lemma $6.1$, $\text{sp}(\psi_{r, \delta}) = \text{supp} \ \hat{\psi}_{r, \delta} \subseteq \mathbb{T} \setminus (-\frac{\delta}{2}, \frac{\delta}{2})$. By definition of $\psi_{r, \delta}$, we have
\begin{align}
\psi_{r, \delta}(n) &= \int_{\mathbb{T}} \hat{\psi}_{r, \delta}(\theta)e^{-in\theta}\frac{d\theta}{2\pi} \notag \\
&= \int_{\mathbb{T}} \hat{\psi}_r(\theta)e^{-in\theta}\frac{d\theta}{2\pi} - \int_{\mathbb{T}} \hat{\psi}_r(\theta) \hat{C}_{\delta}(\theta)e^{-in\theta}\frac{d\theta}{2\pi} \notag \\
&= \psi_r(n) - (\psi_r * C_{\delta})(n). \notag
\end{align}
Letting $r \to 1^-$, we obtain
\[\psi_{\delta}(n) := \lim_{r \to 1^-} \{\psi_r(n) - (\psi_r * C_{\delta})(n)\} = \psi(n) - (\psi * C_{\delta})(n) \ (n \in \mathbb{Z})\]
and 
\[\text{sp}(\psi_{\delta}) \cap \left(-\frac{\delta}{2}, \frac{\delta}{2}\right) = \emptyset \]
by Lemma 4.4. For each $r \in (0, 1)$ and $n \in \mathbb{Z}$, we have
\begin{align}
|\psi_r * C_{\delta}(n)| &= \left|\int_{\mathbb{T}} \hat{\psi}_r \hat{C}_{\delta}(\theta)e^{-in\theta}\frac{d\theta}{2\pi}\right| \notag \\
&\le \int_{\mathbb{T}} |\hat{\psi}_r(\theta)\hat{C}_{\delta}(\theta)|\frac{d\theta}{2\pi} \notag \\
&\le \int_{-\delta}^{\delta} |\hat{\psi}_r(\theta)| \frac{d\theta}{2\pi} < \varepsilon. \notag 
\end{align}
Taking limit as $r \to 1^-$, we obtain $|\psi * C_{\delta}(n)| \le \varepsilon$. Then, it follows that
\[|\psi_{\delta}(n) - \psi(n)| = |\psi * C_{\delta}(n)| < \varepsilon \ (n \in \mathbb{Z}), \]
which means $\|\psi_{\delta}-\psi\|_{\infty} \le \varepsilon$. Therefore, by Theorem $6.1$, we conclude that $\psi \xrightarrow{ac} 0$. We complete the proof.
\end{prf}

\begin{thm}
Let $\psi$ be in $L^{\infty}(\mathbb{Z})$ which vanishes on the negative integers and denote $a_n = \psi(n)$ for $n \ge 0$. Let $\hat{\psi}(z)$ be the function on the unit disc $\mathbb{D}$ defined by
\[\hat{\psi}(z) = \sum_{n=0}^{\infty} a_nz^n. \]
If $\hat{\psi}(z) - \frac{\alpha}{1-z}$ has $H^1$ boundary behavior at $z = 1$, then $\psi$ is one-sidedly almost convergent to $\alpha$, that is, 
\[\lim_{k \to \infty} \frac{1}{k} \sum_{i=0}^{k-1} \psi(n+i) = \alpha \]
uniformly in $n \ge 0$.
\end{thm}

\begin{prf}
Let $\psi_1(n) = \psi(n) - \alpha \cdot I_{\mathbb{Z}_+}(n)$. Then, we have $\hat{\psi}_1(z) = \hat{\psi}(z) - \frac{\alpha}{1-z}$. Hence, by Theorem $7.1$, we obtain $\psi_1 \xrightarrow{ac} 0$. Since, $\psi_1(n) = 0$ for $n < 0$, we have $\psi_1 \xrightarrow{oac} 0$ by Lemma 7.1. Then, by Theorem 7.2 we obtain
\[\frac{1}{k} \sum_{i=0}^{k-1} \psi(n+i) = \alpha \]
uniformly in $n \ge 0$. We complete the proof.
\end{prf}

\begin{cor}
Let $\psi$ be in $l_{\infty}$ and denote $a_n = \psi(n)$ for $n \ge 0$. Let $\hat{\psi}(z)$ be the function on the unit disc $\mathbb{D}$ defined by
\[\hat{\psi}(z) = \sum_{n=0}^{\infty} a_nz^n. \]
If $\hat{\psi}$ has a pole of order $1$ at $z=1$, then $\psi \xrightarrow{oac} -\alpha$, where $\alpha$ is the residue of $\hat{\psi}(z)$ at $z=1$.
\end{cor}

\begin{prf}
By the assumption, on a neighborhood $U$ of $1$, we can write as 
\[\hat{f}(z) = \frac{\alpha}{z-1} + g(z), \]
where $g(z)$ is a regular function on $U$. Hence, we obtain
\[\hat{f}(z) - \frac{-\alpha}{1-z} = g(z). \]
Since the right-hand side of the above equation is continuous on $U$, the function $\hat{f}(z) - \frac{-\alpha}{1-z}$ has $H^1$ boundary behavior at $z=1$ . We obtain the result by Theorem $8.2$.
\end{prf}

\begin{rem}
Note that using Theorem 3.1 of \cite{Deb}, we can relax the assumption in the above theorems that the coefficients $\{a_n\}_{n \ge 0}$ is in $l_{\infty}$ to the assumption that $\liminf_{n \to \infty} (b_{n+1}-b_n) \ge -C$ and $\liminf_{n \to \infty} (c_{n+1}-c_n) \ge -C$ for some constant $C > 0$, where $a_n = b_n + ic_n \ (n \ge 0)$.
\end{rem}

The following can be viewed as a generalization of the classical Fatou-Riesz theorem. We say that a sequence $\{a_n\}_{n \ge 0}$ of complex numbers to be bounded from below if there exists a positive constant $C > 0$ such that $b_n \ge -C$ and $c_n \ge -C$ for every $n \ge 0$, where $b_n$ and $c_n$ are real and imaginary parts of $a_n$, respectively.
\begin{thm}
Let $\{a_n\}_{n \ge 0}$ be a sequence of complex numbers which is bounded from below. Let $f(z)$ be the function on the unit disc $\mathbb{D}$ defined by
\[f(z) = \sum_{n=0}^{\infty} a_nz^n. \]
If $f(z)$ is analytic at the point $z = 1$, then we have $s_n \xrightarrow{oac} f(1)$, where $s_n = \sum_{k=0}^n a_k$.
\end{thm}

\begin{prf}
First, notice that boundedness of the sequence $\{s_n\}_{n \ge 0}$ follows from Theorem $6.2 \ \mathrm{(i)}$ of \cite{Deb}. Next, we show $s_n \xrightarrow{oac} f(1)$. Note that
\[g(z) = \frac{f(z)}{1-z} = \sum_{n=0}^{\infty} s_nz^n \]
and since $f(z)$ is analytic at $z=1$, $g(z)$ has a pole of order $1$ at $z=1$ where its residue is $-f(1)$. Thus, by Corollary $8.1$, we obtain $s_n \xrightarrow{oac} f(1)$. We complete the proof.
\end{prf}

As a corollary of Theorem 8.3, we obtain the Fatou-Riesz theorem on the convergence of power series on the circle of convergence. We need the following elementary lemma which is induced immediately by Theorem 3.2.
\begin{lem}
Let $\{a_n\}_{n \ge 0}$ be in $l_{\infty}$ and $\alpha$ be a complex number. Then, $\lim_{n \to \infty} a_n = \alpha$ if and only if $a_n \xrightarrow{oac} \alpha$ and $\lim_{n \to \infty} (a_{n+1} - a_n) = 0$.
\end{lem}

\medskip
\begin{thm}[Fatou, 1906]
Let $\{a_n\}_{n \ge 0}$ be in $l_{\infty}$ such that $\lim_{n \to \infty} a_n = 0$. Let $f(z)$ be the function on the unit disc $\mathbb{D}$ defined by
\[f(z) = \sum_{n=0}^{\infty} a_nz^n. \]
If $f(z)$ is analytic at the point $z = 1$, then we have $\sum_{n=0}^{\infty} a_n = f(1)$.
\end{thm}

\begin{prf}
By Theorem $8.3$, we obtain $s_n \xrightarrow{oac} f(1)$, where $s_n = \sum_{k=0}^n a_k$. Furthermore, by the assumption, we have $s_n - s_{n-1} = a_n \rightarrow 0$ as $n \to \infty$. Combining Lemma $8.1$, we conclude that $\lim_{n \to \infty} s_n = \sum_{n=0}^{\infty} a_n = f(1)$.
\end{prf}

Now we mention an interesting result of Katznelson and Tzafriri, which can be regarded as a dual of Theorem 8.1. (see \cite{Kat}, \cite{Kor}). 
\begin{thm}[Katznelson and Tzafriri, 1986]
Let $\{a_n\}_{n \ge 0}$ be a bounded sequence of complex numbers. Suppose that the analytic function
\[f(z) = \sum_{n=0}^{\infty} a_nz^n \]
has $H^1$ boundary behavior on the unit circle except for $z=1$. Then, the equation
\[\lim_{n \to \infty} (a_{n+1} - a_n) = 0. \]
holds true.
\end{thm}

Analogous results concerning $G = \mathbb{R}$ can be obtained. Since proofs are similar to the case $G = \mathbb{Z}$ above, in what follows, we exhibit only results without their proofs. First, we define $H^1$ boundary behavior for the functions defined on the right-half plane.
\begin{dfn}
We say that a function $f(z)$ defined on the right-half plane $\mathbb{C}^+ := \{z \in \mathbb{C} : Re (z) > 0\}$ has $H^1$ boundary behavior at the point $z_0 = iy_0 \ (y_0 \in \mathbb{R})$ if there exist a number $\delta > 0$ and a function $F(y)$ in $L^1(y_0-\delta, y_0 + \delta)$ such that $f(x+iy)$ converges to $F(y)$ in $L^1(y_0-\delta, y_0 + \delta)$ as $x \to 0^+$.
\end{dfn}

\begin{thm}
Let $\psi$ be in $L^{\infty}(\mathbb{R})$ which vanishes on the negative reals. Let $\mathscr{L}\psi(s)$ be the Laplace transform of $\psi$, the analytic function on the right half plane $\mathbb{C}^+$ defined by
\[\mathscr{L}\psi(s) = \int_{\mathbb{R}_+} \psi(t)e^{-st}dt, \quad s \in \mathbb{C}^+. \]
If $\mathscr{L}\psi(s)$ has $H^1$ boundary behaviour at $s = 0$, then $\psi$ is almost convergent to $0$.
\end{thm}

\begin{cor}
Let $\psi$ be in $L^{\infty}(\mathbb{R}_+)$ and let $\mathscr{L}\psi(s)$ be its Laplace transform.
If $\mathscr{L}\psi(s) - \frac{\alpha}{s}$ has $H^1$ boundary behaviour at $s = 0$, then $\psi$ is one-sidedly almost convergent to $\alpha$, that is,  
\[\lim_{\theta \to \infty} \frac{1}{\theta} \int_x^{x+\theta} \psi(t)dt = \alpha \]
uniformly in $x \ge 0$.
\end{cor}

\begin{cor}
Let $\psi$ be in $L^{\infty}(\mathbb{R}_+)$ and let $\mathscr{L}\psi(s)$ be its Laplace transform. If $\mathscr{L}\psi(s)$ has a pole of order $1$ at $s = 0$, then $\psi \xrightarrow{oac} \alpha$, where $\alpha$ is the residue of $\mathscr{L}\psi(s)$ at $s=0$.
\end{cor}

\begin{rem}
Using Theorem 3.1 of \cite{Deb} again, we can relax the condition in the above theorems that $\psi$ is in $L^{\infty}(\mathbb{R}_+)$ to the condition that $\psi$ is boundedly decreasing; that is, there are some constants $C > 0$, $\delta > 0$ and $x_0 > 0$ such that $Re \; \psi(x+h) - Re \; \psi(x) \ge -C$ and $Im \; \psi(x+h) - Im \; \psi(x) \ge -C$ whenever $x \ge x_0$ and $h \in [0, \delta]$.
\end{rem}

The following results can be viewed as counterparts of Theorems 8.3 and 8.4 in the continuous setting. Here, we say that $\psi(x)$ is bounded from below if there exists a positive constant $C > 0$ such that $Re \; \psi(x) \ge -C$ and $Im \; \psi(x) \ge -C$ for every $x \ge 0$. We note that boundedness of the function $\Psi(x)$ below can be deduced from Theorem 3.1 of \cite{Deb}.
\begin{thm}
Let $\psi$ be a locally integrable function on $\mathbb{R}_+$ which is bounded from below. Let $\mathscr{L}\psi(s)$ be its Laplace transform. If $\mathscr{L}\psi(s)$ converges in the right half plane and is analytic at $s=0$, then $\Psi \xrightarrow{oac} \mathscr{L}\psi(0)$, where
\[\Psi(x) = \int_0^x \psi(t)dt \ (x \ge 0). \]
\end{thm}

\begin{cor}
Let $\psi$ be a locally integrable function on $\mathbb{R}_+$ such that $\lim_{x \to \infty} \psi(x) = 0$ and let $\mathscr{L}\psi(s)$ be its Laplace transform. If $\mathscr{L}\psi(s)$ is analytic at $s=0$, then $\lim_{x \to \infty} \Psi(x) = \mathscr{L}\psi(0)$, where the function $\Psi$ is defined as above.
\end{cor}

\medskip
Finally, we remark on the relation among Tauberian theorems concerning the behavior of an analytic function $f(z) = \sum_{n=0}^{\infty} a_nz^n$ or that of $\mathscr{L}\psi(s)$ near its boundaries. First, we mention the following classical and famous result (See \cite{Kor}).
\begin{thm}[Hardy-Littlewood]
Let $\{a_n\}_{n \ge 0}$ be a sequence of complex numbers which is bounded from below. Suppose that $f(x) = \sum_{n=0}^{\infty} a_nx^n$ converges for $|x| < 1 \ (x \in \mathbb{R})$ and 
\[\lim_{x \to 1^-} (1-x)f(x) = \alpha. \]
Then, the following formula holds true.
\[\lim_{k \to \infty} \frac{1}{k}\sum_{i=0}^{k-1} a_i = \alpha. \]
\end{thm}
Let $H^1(\mathbb{D})$ be the Hardy space on the unit disc and $H^1_{\{1\}}(\mathbb{D})$ be the class of functions whose elements are the analytic functions on $\mathbb{D}$ having $H^1$ boundary behavour at the point $z=1$. Then, we can summarize as follows: Here let $\{a_n\}_{n \ge 0}$ be a bounded sequence of complex numbers and $f(z) = \sum_{n=0}^{\infty} a_nz^n$.
\[
\begin{cases}
f(x) - \frac{\alpha}{1-x} = o(\frac{1}{1-x}) \ (x \rightarrow 1^-) &\Longrightarrow \lim_{k \to \infty} \frac{1}{k}\sum_{i=0}^{k-1} a_i = \alpha, \\
f(z) - \frac{\alpha}{1-z} \in H^1_{\{1\}}(\mathbb{D}) &\Longrightarrow \lim_{k \to \infty} \frac{1}{k} \sum_{i=0}^{k-1} a_{n+i} = \alpha \ \text{uniformly in $n \in \mathbb{Z}_+$}, \\
f(z) - \frac{a}{1-z} \in H^1(\mathbb{D}) &\Longrightarrow \lim_{k \to \infty} a_k = \alpha. 
\end{cases}
\]
The first one follows from Theorem $8.8$. The second is due to Theorem $8.2$. The third is by the Riemann-Lebesgue lemma. Hence, considering Lemma 8.1, we see that the third assertion above is decomposed into the two assertions of Theorems 8.2 and 8.5.

The corresponding results for the case of $\mathbb{R}$ is in order: First, the right half plane version of Hardy-Littlewood's theorem reads as follows (see \cite{Kor}):
\begin{thm}
Let $\psi$ be a locally integrable function on the half line $\mathbb{R}_+$ which is bounded from below. If $\mathscr{L}\psi(x)$ exists for all $x > 0$ and 
\[\lim_{x \to 0^+} x\mathscr{L}\psi(x) = \lim_{x \to 0^+} x\int_0^{\infty} e^{-xt}\psi(t)dt = \alpha, \]
then we have
\[\lim_{\theta \to \infty} \frac{1}{\theta} \int_0^{\theta} \psi(t)dt = \alpha. \]

\end{thm}

Let $H_{\{0\}}^1(\mathbb{C}^+)$ and $H^1_{\mathbb{R}}(\mathbb{C}^+)$ be the set of analytic functions on the right half plane $\mathbb{C}^+$ having $H^1$ boundary behavior at $s=0$ and the whole line $\mathbb{R}$, respectively. Then, the following result holds true. Here, $\psi$ is in $L^{\infty}(\mathbb{R}_+)$.
\[
\begin{cases}
\mathscr{L}\psi(x) - \frac{\alpha}{x} = o(\frac{1}{x}) \ (x \rightarrow 0^+) &\Longrightarrow \lim_{\theta \to \infty} \frac{1}{\theta} \int_0^{\theta} \psi(t)dt = \alpha, \\
\mathscr{L}\psi(s) - \frac{\alpha}{s} \in H_{\{0\}}^1(\mathbb{C}^+) &\Longrightarrow \lim_{\theta \to \infty} \frac{1}{\theta} \int_x^{x+\theta} \psi(t)dt = \alpha \ \text{uniformly in $x \in \mathbb{R}_+$}, \\
\mathscr{L}\psi(s) - \frac{\alpha}{s} \in H_{\mathbb{R}}^1(\mathbb{C}^+) &\Longrightarrow w^*\mathchar`-\lim_{x \to \infty} \psi_x = \alpha.
\end{cases}
\]
We remark that the third one follows from the Wiener-Ikehara theorem (see for example \cite{Ing}, \cite{Kor}). Briefly, we can say that Theorem 8.6 is an intermediate result between Hardy-Littlewood's theorem and the Wiener-Ikehara theorem.

\end{document}